\documentclass[%
a4paper,%
]%
{article}

\usepackage{amsmath}
\usepackage{amsthm}

\usepackage[matrix,arrow]{xy}

\usepackage{pxfonts}

\usepackage{xspace}
\newcommand{\N}{\ensuremath {\mathbb{N}}\xspace}
\newcommand{\R}{\ensuremath {\mathbb{R}}\xspace}
\newcommand{\Z}{\ensuremath {\mathbb{Z}}\xspace}
\newcommand{\Ci}{\ensuremath {C^\infty}\xspace}
\DeclareMathOperator{\map}{Map}
\newtheoremstyle{thm}{3pt}{3pt}{\itshape}{}{\bfseries}{}{.5em}{}
\newtheoremstyle{thmsub}{3pt}{3pt}{\upshape}{}{\bfseries}{}{.5em}{}
\theoremstyle{thm}
\newtheorem{theorem}{Theorem}[section]
\newtheorem{lemma}[theorem]{Lemma}
\newtheorem{proposition}[theorem]{Proposition}
\newtheorem{defn}[theorem]{Definition}
\newtheorem{corollary}[theorem]{Corollary}
\newtheorem{thm}{Theorem}
\newtheorem{cor}[thm]{Corollary}
\theoremstyle{thmsub}

\newcommand{\m}[1]{\ensuremath {\mathcal{#1}}\xspace}
\newcommand{\hyp}{-\hspace{0pt}}

\newcommand{\abs}[1]{\left\lvert #1 \right\rvert}
\newcommand{\norm}[1][\cdot]{\left\lVert #1 \right\rVert}

\numberwithin{equation}{section}

\DeclareMathOperator{\homeo}{Homeo}
\DeclareMathOperator{\diff}{Diff}
\newcommand{\ci}{c^\infty}

\begin{document}

\title{%
Constructing Smooth Loop Spaces%
}
\date{\today}
\author{Andrew Stacey}

\maketitle

\begin{abstract}
  We consider the general problem of constructing the structure of a
  smooth manifold on a given space of loops in a smooth finite
  dimensional manifold.  By generalising the standard construction for
  smooth loops, we derive a list of conditions for the model space
  which, if satisfied, mean that a smooth structure exists.

  We also show how various desired properties can be derived from the
  model space; for example, topological properties such as
  paracompactness.  We pay particular attention to the fact that the
  loop spaces that can be defined in this way are all homotopy
  equivalent; and also to the action of the circle by
  rigid rotations.
\end{abstract}


\section{Introduction}
\label{sec:intro}

It is often convenient to regard a space of certain loops in a smooth
manifold as a smooth manifold itself with the aim of doing
differential topology thereon.  Depending on the application this
approach can vary from the conceptual to the rigorous.  The two most
popular types of loop are continuous and smooth, for both of which
there is a rigorous theory of infinite dimensional manifolds making
these into smooth manifolds: \cite{wk}, \cite{sl}, \cite{jm},
\cite{ho}.  Other types of loop have also been considered: it is often
convenient to use a manifold modelled on a Hilbert space whence one
usually uses the space of loops with square-integrable first
derivative.

There is a standard method of constructing the smooth structure which
is used in each case mentioned above.  Our main purpose in this paper
is to extend this construction to a reasonably arbitrary class of
loops.  In so doing, we obtain a list of conditions on the model space
such that if they hold then this method applies.  This enables us to
reduce the general problem of whether or not a particular type of loop
forms a smooth manifold to a check-list for the model space and means
that we can avoid writing out the full construction in each and every
case.

Before giving the list of conditions, we feel it relevant to comment
on calculus.  The examples of spaces given in the first paragraph were
all ``nice'' as regards calculus.  Two were Banach spaces (one a
Hilbert space, no less) and the other is one of the nicest Fr\'echet
spaces that one could hope to meet.  In all of these cases, calculus
is well-understood and well-defined.  However, once one leaves the
realm of Fr\'echet spaces and departs for more general locally convex
topological spaces then the notion of what is ``smooth'' becomes
increasingly hard to pin down.  The usual idea of taking ``smooth'' to
mean ``infinitely differentiable'' leads to many complications, not
least that this is not uniquely defineable.  Fortunately, an
alternative approach has been developed in which the notion of
``smooth'' is based on something other than differentiability.  This
calculus is laid out in the weighty tome \cite{akpm}.  The
introduction and the historical remarks at the end of chapter~1 of
\cite{akpm} make for an interesting read on this subject.

The impact that this has on our work is more subtle than might be
expected.  The place where one would expect this issue to arise is in
showing that the transition functions are smooth.  However, this
depends on certain functions on the model space being smooth and so we
build this into one of our conditions.  We can therefore phrase the
corresponding condition in such terms that it could apply whatever
type of calculus we were using.

There are other places, however, where the calculus used makes an
appearance.  The most important being the definition of an infinite
dimensional manifold.  One of the issues in infinite dimensional
calculus is that maps can be smooth without being continuous and this
leads to a certain \emph{laissez faire} attitude to topology.  The
traditional definition of a manifold is of a topological space with a
smooth atlas.  The definition in \cite{akpm} is of a \emph{set} with a
smooth atlas which is then topologised using said atlas.  Thus if one
wishes to apply the results of this paper using a calculus other than
that of \cite{akpm}, this issue might be important.  Certainly, the
traditional approach to building infinite dimensional manifolds
modelled on Banach spaces has been to mirror the standard approach of
topology first and smooth structure second.  Therefore, if taking a
different calculus, it might be necessary to add the additional step
of checking that the original topology and the new topology were one
and the same.  We shall not bother with this issue explicitly, but
shall provide some tools which will help if it is considered important
by others. 

Having made that point, we turn to our list of conditions.  We start
with a class of maps \(S^1 \to \R\) which we write as \(L^x \R\), or
as \(C^x(S^1, \R)\) when we want to emphasise the domain, and refer to
as \(C^x\)\hyp{}loops.  We want to consider actual maps, and not
equivalence classes of maps, because we want to be able to ``locate''
our maps on a manifold in order to be able to examine them in charts.
Thus we regard \(L^x \R\) as a subset of \(\map(S^1,\R)\).  Using the
natural identification of \(\map(S^1, \R^n)\) with \(\map(S^1, \R)^n\)
we define \(L^x \R^n\) as \((L^x \R)^n\).  For a subset \(A \subseteq
\R^n\), we define \(L^x A\) as the subset of \(L^x \R^n\) consisting
of maps which take values in \(A\).

Our conditions are:
\begin{enumerate}
\item Being in \(L^x \R\) is a \emph{local} property.
\label{cond:local}

That is, for a loop \(\gamma \colon S^1 \to \R\), then \(\gamma \in
L^x \R\) if there exists an open cover \(\m{U}\) of \(S^1\) and
loops \(\gamma_U \in L^x \R\) for \(U \in \m{U}\) such that
\(\gamma\) agrees with \(\gamma_U\) on \(U\).

\item The set \(L^x \R\) is a subspace of \(\map(S^1, \R)\).
\label{cond:vspace}

\item The vector space \(L^x \R\) can be given a topology with
respect to which it is a locally convex topological vector space.
\label{cond:lctvs}

\item The locally convex topological vector space \(L^x \R\) is
\emph{convenient}.
\label{cond:cmplt}

This is a completeness condition.  We have phrased it in the language
of \cite{akpm} but it is the same as a concept known as \emph{locally
  complete} which is from ordinary functional analysis.  Local
completeness is weaker than sequential completeness, though it
coincides with completeness for metrisable spaces.  This completeness
condition is to ensure that derivatives that ought to exist actually
do.

\item As subspaces of \(\map(S^1, \R)\) we have inclusions:
\[
  L \R \subseteq L^x \R \subseteq L^0 \R
\]
where \(L \R = \Ci(S^1, \R)\) and \(L^0 \R = C^0(S^1, \R)\).  The
inclusion maps are all continuous when each is given its standard
topology.
\label{cond:smthcts}

We will be considering loops in a smooth manifold and therefore will
need to know that the condition of being a \(C^x\)--loop is invariant
under post-composition by diffeomorphisms.  This essentially forces
smooth loops to be \(C^x\)--loops.  For the other inclusion, as we
remarked above we want to be able to ``locate'' a loop on a manifold
so that we can consider it in charts.  The simplest way to do this is
to ensure that a \(C^x\)--loop is continuous.

\item The action of post-composition of \(C^x\)--loops by smooth maps
is well-defined and is smooth.  That is, let \(U \subseteq \R^m\) and
\(V \subseteq \R^n\) be open sets; let \(\phi \colon U \to V\) be a smooth
map.  The induced map \(\phi_* \colon L^x U \to L^x V\),
\(\gamma \mapsto \phi \circ \gamma\), is well-defined and is smooth.
\label{cond:postcomp}

This is the crucial condition that will ensure that the transition
functions are defined and are diffeomorphisms.
\end{enumerate}

Having stated our conditions, we can now state our first theorem.  We
make two assumptions on our target manifold: that it be orientable and
that it have no boundary.  The first is really a convenience to allow
us not to have to discuss twisted model spaces; the second is
necessary as the loop space of a manifold with boundary is a
complicated object indeed.

\begin{thm}
\label{th:smooth}
Let \(L^x \R\) be a class of maps satisfying the above
conditions.  Let \(M\) be a smooth, orientable finite dimensional
manifold without boundary.  Then \(L^x M\) can be defined and is
a smooth manifold in the sense of \cite{akpm}.
\end{thm}

We emphasise that the phrase ``in the sense of \cite{akpm}'' does not
refer to the calculus but to the definition of a smooth manifold once
one has decided on a calculus.

Having defined the smooth structure, the next question is to examine
the general properties of the manifold.  In most cases, these descend
from the model space.  The result which allows us to devolve these
properties is the following theorem on submanifolds.

\begin{thm}
\label{th:submfd}
Let \(L^x \R\) be a class of maps satisfiying the above
conditions.  Let \(M, N\) be smooth, orientable finite dimensional
manifolds without boundary and suppose that there is an embedding of
\(M\) as a submanifold of \(N\).  Then \(L^x M\) 
is an embedded submanifold of \(L^x N\).  Moreover, if \(M\) is
closed, resp.\ open, in \(N\) then \(L^x M\) is closed, resp.\
open, in \(L^x N\).
\end{thm}

\begin{cor}
\label{cor:top}
In the statement of theorem~\ref{th:submfd} suppose that we can take
\(N = \R^n\) with \(M\) closed in \(N\).  Then the following
properties are inherited by \(L^x M\) from \(L^x \R\):
separable, metrisable, Lindel\"of, paracompact, normal, smoothly
regular, smoothly paracompact, and smoothly normal.

That is, those properties which hold for \(L^x \R\) also hold
for \(L^x M\).
\end{cor}

The last property that we wish to examine is the natural circle
action, and more generally the natural action of the diffeomorphisms
of the circle.

\begin{thm}
\label{th:diffact}
Under the conditions of corollary~\ref{cor:top}, the actions of the
circle and of the diffeomorphisms of the circle are also inherited by
\(L^x M\) from \(L^x \R\).
\end{thm}

In light of this, we conclude this paper with a discussion as to the
various possible levels of continuity and smoothness of the circle
acting on an infinite dimensional locally convex topological vector
space (lctvs).

\medskip

This paper is organised as follows.  In section~\ref{sec:prelim} we
prove some preliminary results, including the definition of \(L^x
M\).  Section~\ref{sec:smooth} is concerned with the construction of
the charts and showing that the transition maps are diffeomorphisms;
in particular we prove theorem~\ref{th:smooth}.
In section~\ref{sec:topology} we transfer our attention to the
topology of the manifold and prove theorem~\ref{th:submfd} and its
corollaries.  Finally, in section~\ref{sec:diff} we look at the action
of the circle and its diffeomorphisms.

\medskip

The standard construction of the smooth structure on the space of
smooth loops can be found in many places, for example in \cite{pm3}
and in \cite{akpm}.  Some other articles and books which deal with the
infinite dimensional manifolds in varying levels of generality are:
\cite{jm}, \cite{ho}, \cite{je3}, \cite{jeke}, \cite{sl}, \cite{wk}.
Most of the work in this paper is firmly in the realms of differential
topology and should be comprehensible to anyone with a firm grasp of
the basics of the theory in finite dimensions.  The exception to this
is the discussion of actions of the circle and diffeomorphism group
which uses some standard functional analysis.  This may be unfamiliar
to differential topologists, at whom this article is aimed, in which
case we recommend \cite{hs} and \cite{hj} for the necessary
background.

\medskip

We regard the circle as the quotient \(\R / \Z\) and shall write it
additively.  We shall often write a small neighbourhood of a point as
\((t - \epsilon, t + \epsilon)\) without worrying about
``wrap-around'';  either the ``wrap-around'' will have no effect on
the subsequent discussion or we will be allowed to take \(\epsilon\)
small enough that there is no ``wrap-around''.  We shall also employ
the language of \emph{intervals} for connected subsets of \(S^1\) --
including \(S^1\) itself.

\section{Preliminaries}
\label{sec:prelim}

In this section we shall set up the basic machinery that we need to
define and construct the smooth manifold of loops.  From hereon, let
us assume that we have a class of loops, \(L^x \R\), satisfying
the conditions stated in the introduction.  Let \(M\) be a smooth,
orientable, finite dimensional manifold without boundary.  Our first
task is to define \(L^x M\).  Our second is to define and
examine the space of \(C^x\)\hyp{}sections in a smooth vector bundle
over \(S^1\); these will prove crucial in the atlas for \(L^x
M\).

\subsection{Loops in a Manifold}
\label{sec:defloops}

To define \(L^x M\) we need to show that we can examine a loop
locally to see whether or not it is in \(L^x M\).
Condition~\ref{cond:local} is almost what we need but isn't quite
local enough.

\begin{defn}
Let \(I \subseteq S^1\) be an open interval.  Define \(C^x(I, \R)\) to
be the space of maps \(\gamma \colon I \to \R\) which are locally of type
\(C^x\).  That is, there is an open cover \(\m{U}\) of \(I\) and maps
\(\gamma_U \in L^x \R\) for \(U \in \m{U}\) such that \(\gamma\)
agrees with \(\gamma_U\) on \(U\).
\end{defn}

Note that we don't assume that the whole function extends, merely that it
locally extends.  It follows from the definition that the restriction
map \(C^x(I, \R) \to C^x(J, \R)\) is defined for \(J \subseteq I\).

\begin{lemma}
\label{lem:local}
In the locality condition, it is enough to assume that the local
functions are only defined locally.

That is, a map \(\gamma \colon S^1 \to \R\) is a \(C^x\)\hyp{}map if there
is a cover \(\m{U}\) of \(S^1\) by open intervals and functions
\(\gamma_U \in C^x(U, \R)\) such that \(\gamma\) agrees with
\(\gamma_U\) on \(U\).
\end{lemma}

\begin{proof}
Let \(t \in S^1\).  Then there is some \(U \in \m{U}\) with \(t \in
U\).  As \(\gamma_U \in C^x(U,\R)\) there is an open cover \(\m{V}\)
of \(U\) and loops \(\beta_V \in L^x \R\) such that \(\gamma_U\)
agrees with \(\beta_V\) on \(V\).  There is some \(V \in \m{V}\) such
that \(t \in V\).  Then on \(V\), \(\beta_V\) agrees with \(\gamma_U\)
which agrees with \(\gamma\).  Repeating this for all \(t \in S^1\)
gives the family required to apply condition~\ref{cond:local}.
\end{proof}

Another piece of preliminary work that we need to do, or rather just
to note as it is trivial, is to show that all of our conditions and
results are equally valid for \(\R^n\) as for \R.
Condition~\ref{cond:postcomp} is already in full generality.

\begin{lemma}
Let \(L^x \R\) be a class of maps satisfying the conditions of
section~\ref{sec:intro}.  Then \(L^x \R^n\) satisfies analogous
conditions and the corresponding version of lemma~\ref{lem:local}.
\end{lemma}

\begin{proof}
This is trivial and follows from the fact that \(L^x \R^n\) is
canonically identified with \((L^x \R)^n\).
\end{proof}

One other result that we need, which is equally trivial, is the
following statement about open sets.

\begin{lemma}
Let \(U \subseteq \R^n\) be open.  Then \(L^x U\) is open in
\(L^x \R^n\).
\end{lemma}

\begin{proof}
This holds for \(L^0 \R^n\) and so holds because \(L^x U =
L^0  U \cap L^x \R^n\).
\end{proof}

With this in place we can define our space of interest.

\begin{defn}
Let \(L^x \R\) be a class of maps satisfying the conditions of
section~\ref{sec:intro}.  Let \(M\) be a smooth finite dimensional
manifold.  Define \(L^x M\) to be the subset of \(\map(S^1, M)\)
consisting of those loops \(\gamma \colon S^1 \to M\) for which there
exists a covering \(\{I_\alpha \colon \alpha \in A\}\) of \(S^1\) by open
intervals and charts \(\{(\iota_\alpha, U_\alpha, V_\alpha) : \alpha
  \in A\}\) for \(M\) (not necessarily making a full atlas) such that
for each \(\alpha \in A\), \(\gamma(I_\alpha) \subseteq V_\alpha\) and
the map:
\[
  \gamma_\alpha \colon I_\alpha \to S^1 \xrightarrow{\gamma} V_\alpha
  \xrightarrow{{\iota_{\alpha}}^{-1}} U_\alpha
\]
lies in \(C^x(I_\alpha, U_\alpha)\).
\end{defn}

Thus we have defined \(C^x\)--loops in \(M\) to be those that look
like \(C^x\)--loops whenever we look locally.  Lemma~\ref{lem:local}
and condition~\ref{cond:postcomp} easily combine to show that this
definition does not depend on any of the choices made.

There is another way to make this definition; if \(M\) were a
submanifold of, say, \(\R^n\) then we already have a definition of
\(L^x M\): namely that subset of \(L^x \R^n\) of loops
which take values in \(M\).  The next result shows that these two
definitions coincide.  We prefer the above as the actual definition as
it treats the manifold without reference to any surrounding space.

\begin{proposition}
\label{prop:cxsub}
Let \(L^x \R\) be a class of maps satisfying the conditions in
section~\ref{sec:intro}.  Let \(M,N\) be smooth, finite dimensional
manifolds with \(M\) an embedded submanifold of \(N\).  Then 
\[
  L^x M = \{\gamma \in L^x N : \gamma(S^1) \subseteq M\}.
\]
\end{proposition}

\begin{proof}
Since this is true for arbitrary maps, what we need to show is that a
loop in \(M\) is a \(C^x\)\hyp{}loop when viewed in \(M\) if and only
if it is a \(C^x\)\hyp{}loop when viewed in \(N\).  To do this, we
ensure that the charts in \(M\) in which we are looking are all
submanifold charts; that is, restrictions of charts on \(N\) which
take \(M\) to \(\R^k\) inside \(\R^n\).  The desired result then
follows from the fact that a loop in \(\R^n\) is a  \(C^x\)\hyp{}loop
if and only if its co-ordinates are \(C^x\)\hyp{}loops; whereupon the
\(C^x\)\hyp{}loops in \(\R^k\) are precisely the \(C^x\)\hyp{}loops in
\(\R^n\) which happen to lie in \(\R^k\).
\end{proof}

\subsection{Sections and Submanifolds}
\label{sec:cxsect}

In this section we define and examine the space of
\(C^x\)\hyp{}sections of a vector bundle over \(S^1\).

\begin{defn}
Let \(E \to S^1\) be a smooth fibre bundle.  Define
\(\Gamma_{S^1}^x(E)\) as the space of sections of \(E\) which are
\(C^x\)\hyp{}loops when viewed as maps into the total space of \(E\).
\end{defn}

In the particular case that \(E\) is an orientable vector bundle, we
can identify this space of sections with \(L^x \R^n\).

\begin{lemma}
\label{lem:cxsec}
Let \(E \to S^1\) be a smooth orientable vector bundle of fibre
dimension \(n\).  A smooth trivialisation of \(E\) defines a bijection
\(\Gamma^x_{S^1}(E) \to L^x \R^n\).  The map \(L^x \R^n
\to L^x \R^n\) induced by two trivialisations of \(E\) is a
linear diffeomorphism.
\end{lemma}

\begin{proof}
As it is obvious that a smooth trivialisation of \(E\) defines a
bijection from the space of all sections of \(E\) to \(\map(S^1,
\R^n)\) all we need to check to show the first part is that a section
of \(E\) is a \(C^x\)\hyp{}section if and only if this bijection takes
it to a \(C^x\)\hyp{}loop in \(\R^n\).

Condition~\ref{cond:postcomp} assures us that a diffeomorphism on the
target space induces a bijection on the spaces of \(C^x\)\hyp{}loops.
Therefore we have a bijection from \(L^x E\) to \(L^x (S^1
\times \R^n)\).  As the trivialisation of \(E\) intertwines the
projection maps, this bijection takes sections to sections and so
induces a bijection \(\Gamma_{S^1}^x(E) \to \Gamma_{S^1}^x(S^1 \times
\R^n)\).  Thus the problem is reduced to the case of a trivial vector
bundle.

Now it is clear from the definition of \(C^x\)\hyp{}loops in a
manifold that a loop in a (finite) product is a \(C^x\)\hyp{}loop if
and only if each of the factors is a \(C^x\)\hyp{}loop.  Therefore a
loop in \(S^1 \times \R^n\) is a \(C^x\)\hyp{}loop if and only if the
projections to \(S^1\) and to \(\R^n\) are \(C^x\)\hyp{}loops.  Now a
loop is a section if and only if it projects to the identity on
\(S^1\) and the identity is smooth, whence a \(C^x\)\hyp{}loop.
Therefore a section of \(S^1 \times \R^n\) is a \(C^x\)\hyp{}section
if and only if the projection to \(\R^n\) produces a
\(C^x\)\hyp{}loop.

Tracing this through shows that the trivialisation does induce a
bijection \(\Gamma_{S^1}^x(E) \to L^x \R^n\).

Two such trivialisations of \(E\) define a diffeomorphism \(\phi
\colon S^1
\times \R^n \to S^1 \times \R^n\) covering the identity on \(S^1\).
We extend this to a smooth map \(\R^2 \times \R^n \to \R^2 \times
\R^n\), viewing \(S^1\) as a submanifold of \(\R^2\).  Note that we do
not assume that this extension is a diffeomorphism (an extension to a
diffeomorphism may not exist).  The induced map \(L^x \R^n \to L^x
\R^n\) factors as:
\begin{align*}
L^x\R^n &\to L^x(\R^2 \times \R^n)  & \gamma &\mapsto
(0,\gamma) \\
L^x(\R^2 \times \R^n) &\to L^x(\R^2
\times \R^n) & (\alpha, \gamma) &\mapsto (\alpha + \iota, \gamma) \\
L^x(\R^2
\times \R^n) &\xrightarrow{\phi_*} L^x(\R^2 \times \R^n)  &
(\alpha, \gamma) &\mapsto \phi \circ (\alpha, \gamma) \\
L^x( \R^2 \times \R^n) &\to
L^x \R^n & (\alpha, \gamma) &\mapsto \gamma,
\end{align*}
where \(\iota \colon S^1 \to \R^2\) is the inclusion.  The first map is
continuous and linear, hence smooth.  The second map is a translation,
hence smooth.  The third map is smooth by
condition~\ref{cond:postcomp}.  The final map is continuous and
linear, hence smooth.  Thus \(\phi_* \colon L^x\R^n \to L^x \R^n\) is
smooth.  Applying the same to \(\phi^{-1}\) shows that \(\phi_*\) is a
diffeomorphism, as required.
\end{proof}

Using this we transfer the smooth structure of \(L^x \R^n\) to
\(\Gamma_{S^1}^x(E)\).  If we are using the calculus of \cite{akpm},
it is a curious fact that although \(\phi_*\) is a linear
diffeomorphism, it need not be a homeomorphism as it, or its inverse,
need not be continuous.  They will, however, be \emph{bounded} maps.
If we assume that the topology on \(L^x \R^n\) is \emph{bornological}
-- a condition that we can readily impose by a standard alteration of
the topology -- then bounded maps are continuous and so we do have a
homeomorphism.

\begin{corollary}
Let \(E \to S^1\) be a finite dimensional orientable smooth vector
bundle.  Then \(\Gamma^x_{S^1}(E)\) naturally has the structure of a
convenient vector space and is diffeomorphic to \(L^x \R^n\), where
\(n = \dim E\). \hspace*{\fill}\qedsymbol
\end{corollary}

We can adapt the proof of lemma~\ref{lem:cxsec} to demonstrate the
following properties of \(\Gamma^x_{S^1}(E)\).

\begin{lemma}
\label{lem:chartprelim}
Let \(E,F \to S^1\) be finite dimensional orientable smooth vector
bundles.  Let \(U \subseteq E\) and \(V \subseteq F\) be open subsets
of the total space and \(\phi \colon U \to V\) a smooth map covering the
identity on \(S^1\).  Let \(\Gamma_{S^1}^x(U) \coloneqq \{\gamma \in
  \Gamma_{S^1}^x(E) : \gamma(S^1) \subseteq U\}\), and similarly for
\(V\).  Then \(\Gamma_{S^1}^x(U)\) is open in \(\Gamma_{S^1}^x(E)\)
and the induced map \(\gamma \mapsto \phi \circ \gamma\) is a smooth
map \(\Gamma_{S^1}^x(U) \to \Gamma_{S^1}^x(V)\).
\end{lemma}

\begin{proof}
It is sufficient to examine the case where \(E\) and \(F\) are
trivial; say, \(E = S^1 \times \R^m\) and \(F = S^1 \times \R^n\).  In
this case we have a topological embedding of \(L^x \R^m\) as an affine
subspace of \(L^x(\R^2 \times \R^m)\) via \(\gamma \mapsto (\iota,
\gamma)\).  There is a set \(W \subseteq \R^2 \times \R^m\) which
restricts to \(U\) on \(S^1 \times \R^m\) and, under the above
embedding, \(\Gamma_{S^1}^x(U)\) is the intersection of \(L^x \R^m\)
with \(L^x W\), hence is open in \(\Gamma_{S^1}^x(E)\).

Now smoothness is a local property so we may assume that \(\phi \colon U
\to V\) extends to a smooth map \(\R^2 \times \R^m \to \R^2 \times
\R^n\).  To deduce the general case from this we choose a sequence of
open sets \(U_n\) such that \(U = \bigcup U_n\) and \(\overline{U}_n
\subseteq U_{n+1}\).  Using bump functions we can define maps \(\phi_n
\colon \R^2 \times \R^m \to \R^2 \times \R^n\) such that \(\phi_n = \phi\)
on \(U_n\).  Thereupon if we can show that each \(\phi_n\) is smooth
then we can deduce that \(\phi\) is locally smooth and hence smooth.

We now use the same method as in the proof of lemma~\ref{lem:cxsec} to
deduce that \(\phi_* \colon \Gamma_{S^1}^x(U) \to \Gamma_{S^1}^x(V)\) is
smooth.
\end{proof}

\section{Building a Smooth Manifold}
\label{sec:smooth}

In this section we construct the charts for \(L^x M\) and show that
the transition maps are smooth.

\subsection{Charts}
\label{sec:charts}

The key tool for defining the charts for the loop space is the notion
of a \emph{local addition} on \(M\), cf~\cite[\S 42.4]{akpm}:
\begin{defn}
\label{def:locadd}
Let \(U \subseteq M\) be an open subset of \(M\).  A \emph{local
  addition over \(U\)} for \(M\) consists of an open subset \(U
\subseteq M\) and smooth map \(\eta \colon T U \to U\) such that
\begin{enumerate}
\item the composition of \(\eta\) with the zero section is the
identity on \(U\), and

\item there exists an open neighbourhood \(V\) of the diagonal in
\(U\) such that the map \(\pi \times \eta \colon T U \to U \times U\) is a
diffeomorphism onto \(V\).
\end{enumerate}

For a subset \(A \subseteq M\), a \emph{local addition for \(A\)} is a
local addition defined over a neighbourhood of \(A\).  If \(f \colon X \to
M\) is a map, a \emph{local addition for \(f\)} is a local addition
defined over a neighbourhood of the image of \(f\).
\end{defn}

This is closely related to what is called a local addition in \cite[\S
  42.4]{akpm} but is not quite the same.  One difference, that we use
the whole of the fibres, is for simplicity whilst the other
difference, that we do not assume it to be defined on the whole of
\(M\), is to make later analysis easier.  The following result is
contained in the discussion following the definition of a local
addition in \cite[\S 42.4]{akpm}:

\begin{proposition}
Any finite dimensional manifold without boundary admits a local
addition defined over the whole of the manifold. \hspace*{\fill}\qedsymbol
\end{proposition}

We start by constructing our chart maps for \(L^x M\).  They
will be anchored at \emph{smooth} loops rather than arbitrary
elements of \(L^x M\).  This is to ensure that all the maps
between finite dimensional objects are smooth so we don't need to
consider \(C^x\)\hyp{}maps with arbitrary domains.

\begin{lemma}
  Let \(\alpha \colon S^1 \to M\) be a \emph{smooth} loop.  Let \(\eta
  \colon T U \to U\) be a local addition for \(\alpha\) with
  neighbourhood \(V\) of the diagonal.  Define the set \(U_\alpha
  \subseteq L^x M\) by:
  \[
  U_\alpha \coloneqq \{\beta \in L^x M \colon (\alpha, \beta) \in L^x V\}.
  \]
  Then \(\pi \times \eta \colon T U \to V\) induces a bijection from
  \(\Gamma^x_{S^1}(\alpha^* T M)\) to \(U_\alpha\).  Under this
  bijection, the zero section of \(\alpha^* T M\) maps to
  \(\alpha\).
\end{lemma}

\begin{proof}
By definition, the image of \(\alpha\) lies in \(U\).  As \(U\) is an
open subset of \(M\), the bundles \(\alpha^* T M\) and \(\alpha^* T
U\) are naturally identified.  We claim that there is a diagram:
\[
\xymatrix{
    L^x T U \ar[r]^{(\pi \times \eta)_*} &
    L^x V \\
  \Gamma^x_{S^1}(\alpha^* T M) \ar[u] &
  U_\alpha \ar[u]_{\beta \mapsto (\alpha, \beta)},
}
\]
such that the map at the top is a bijection and takes the image of the
left-hand vertical map to the image of the right-hand one.  Both of
the vertical maps are injective -- the right-hand one obviously so, we
shall investigate the left-hand one in a moment -- and thus the
bijection \((\pi \times \eta)_*\) induces a bijection from the lower
left to the lower right.

As \(T U\) is an open subset of \(T M\) and \(V\) of \(M \times M\),
both are smooth manifolds.  The map \(\pi \times \eta \colon T U \to V\) is
a diffeomorphism and hence induces a bijection on the sets of
\(C^x\)\hyp{}maps into each.  This is the map we have labelled
\((\pi \times \eta)_*\).

The left-hand vertical map, \(\Gamma^x_{S^1}(\alpha^* T M) \to
L^x T U\), is defined as follows: the total space \(\alpha^* T
M = \alpha^* T U\) is:
\[
\{(t,v) \in S^1 \times T U \colon \alpha(t) = \pi(v)\}.
\]
It is an embedded submanifold of \(S^1 \times T U\).  Therefore by
proposition~\ref{prop:cxsub}, a map into \(\alpha^* T M\) is a
\(C^x\)\hyp{}map if and only if the compositions with the projections
to \(S^1\) and to \(T U\) are both \(C^x\)\hyp{}maps.  Now a map \(S^1
\to \alpha^* T U\) is a section if and only if it projects to the
identity on \(S^1\).  Therefore there is a bijection (of sets):
\begin{align*}
\Gamma^x_{S^1}(\alpha^* T M) &\cong \{\beta \in L^x T U \colon (t,
  \beta(t)) \in \alpha^* T M \text{ for all } t \in S^1\} \\
&= \{\beta \in L^x T U \colon \alpha(t) = \pi \beta(t) \text{ for all
  } t \in S^1\} \\
&= \{\beta \in L^x T U \colon \pi_* \beta = \alpha\}.
\end{align*}
In particular, the map \(\Gamma^x_{S^1}(\alpha^* T M) \to L^x T U\) is
injective.

We apply \((\pi \times \eta)_*\) to the image of
\(\Gamma^x_{S^1}(\alpha^* T M)\) and see that it is the preimage under
this map of everything of the form \((\alpha, \gamma)\) in \(L^x
V\).  By construction, \(\gamma \in L^x M\) is such that
\((\alpha, \gamma) \in L^x V\) if and only if \(\gamma \in
U_\alpha\).  Hence \((\pi \times \eta)_*\) identifies the image of
\(\Gamma^x_{S^1}(\alpha^* T M)\) with \(\{\alpha\} \times U_\alpha\).

Finally, note that the zero section of \(\alpha^* T M\) maps to the
image of \(\alpha\) under the zero section of \(T U\).  Since \(\eta\)
composed with the zero section of \(T U\) is the identity on \(U\),
the image of the zero section of \(\alpha^* T M\) in \(V\) is
\((\alpha,\alpha)\) as required which projects to \(\alpha\) in
\(U_\alpha\).
\end{proof}

\begin{defn}
  Let \(\Psi_\alpha \colon \Gamma^x_{S^1}(\alpha^* T M) \to U_\alpha\) be
  the resulting bijection.
\end{defn}

In detail, this map is as follows: let \(\beta \in
\Gamma^x_{S^1}(\alpha^* T M)\) and let \(\tilde{\beta}\) be the
corresponding loop in \(T U\), so \(\beta(t) = (t,
\tilde{\beta}(t))\).  Then \((\pi \times \eta)_*(\tilde{\beta}) =
(\alpha, \eta_*(\tilde{\beta}))\) so \(\Psi_\alpha(\beta) =
\eta_*(\tilde{\beta})\).

The domains of these charts are naturally convenient vector spaces.
On the other end, we need to show that the codomains cover \(L^x M\).

\begin{lemma}
The codomains of the charts cover \(L^x M\).
\end{lemma}

\begin{proof}
This follows from the density of \(L M\) in \(L^0 M\).
We choose a local addition defined over the whole of \(M\), \(\eta \colon T
M \to M\), with corresponding neighbourhood \(V \subseteq M \times M\)
of the diagonal.  As such, for any \(\beta \in L^0 M\) there is
some \(\alpha \in L M\) such that \((\alpha, \beta) \in
L^0 V\).  Whereupon, if \(\beta \in L^x M\) then \(\beta
\in U_\alpha\).  Hence the sets \(U_\alpha\) cover \(L^x M\).
\end{proof}

\subsection{Transitions}
\label{sec:trans}

Having defined the charts, we turn to the transition functions.  Let
\(\alpha_1, \alpha_2\) be smooth loops in \(M\).  Let \(\eta_1 \colon T U_1
\to U_1\) and \(\eta_2 \colon T U_2 \to U_2\) be local additions for
\(\alpha_1\) and \(\alpha_2\) respectively, with corresponding open
sets \(V_1 \subseteq U_1 \times U_1\) and \(V_2 \subseteq U_2 \times
U_2\).  Let \(\Psi_1 \colon \Gamma^x_{S^1}({\alpha_1}^* T M) \to
U_{\alpha_1}\), \(\Psi_2 \colon \Gamma^x_{S^1}({\alpha_2}^* T M) \to
U_{\alpha_2}\) be the corresponding charts.  Let \(U_{1 2} \coloneqq
U_{\alpha_1} \cap U_{\alpha_2}\).

\begin{lemma}
  Let \(W_{1 2} \subseteq {\alpha_1}^* T M\) be the set:
  \[
  \{(t,v) \in {\alpha_1}^* T M \colon (\alpha_2(t), \eta_1(v)) \in V_2\}.
  \]
  Then \(W_{1 2}\) is open and \({\Psi_1}^{-1}(U_{1 2}) =
  \Gamma^x_{S^1}(W_{1 2})\).
\end{lemma}

Here \(\Gamma^x_{S^1}(W_{1 2})\) is the set of sections of
\({\alpha_1}^* T M\) which take values in \(W_{1 2}\).

\begin{proof}
  The set \(W_{1 2}\) is open as it is the preimage of an open set via
  a continuous map.  To show the second statement we need to prove
  that \(\gamma \in \Gamma^x_{S^1}({\alpha_1}^* T M)\) takes values in
  \(W_{1 2}\) if and only if \(\Psi_1(\gamma) \in U_2\) (by
  construction we already have \(\Psi_1(\gamma) \in U_1\)).

  So let \(\gamma \in \Gamma^x_{S^1}(\alpha^* T M)\) and let
  \(\tilde{\gamma}\) be the image of \(\gamma\) in \(L^x T U\).
  Thus \(\gamma(t) = (t, \tilde{\gamma}(t))\).  Now \(\gamma\) takes
  values in \(W_{1 2}\) if and only if:
  \[
  \big(\alpha_2(t), \eta(\tilde{\gamma}(t))\big) \in V_2
  \]
  for all \(t \in S^1\).  That is to say, if and only if \((\alpha_2,
  \eta_*(\tilde{\gamma})) \in L^x V_2\).  By definition, this
  is equivalent to the statement that \(\eta_*(\tilde{\gamma}) \in
  U_{\alpha_2}\).  Now \(\eta_*(\tilde{\gamma}) = \Psi_1(\gamma)\)
  so \(\gamma\) takes values in \(W_{1 2}\) if and only if
  \(\Psi_1(\gamma) \in U_{\alpha_1} \cap U_{\alpha_2}\).
\end{proof}

\begin{proposition}
  \label{prop:trans}
  The transition function:
  \[
  \Phi_{1 2} \coloneqq {\Psi_1}^{-1} \Psi_2 \colon {\Psi_1}^{-1}(U_{1 2}) \to
  {\Psi_2}^{-1}(U_{1 2})
  \]
  is a diffeomorphism.
\end{proposition}

\begin{proof}
  We define \(W_{2 1} \subseteq {\alpha_2}^* T M\) as the set of
  \((t,v) \in \alpha_2^* T M\) with \((\alpha_1(t), \eta_1(v)) \in
  V_1\).  As for \(W_{1 2}\), \(\Psi_2^{-1}(U_{1 2}) =
  \Gamma^x_{S^1}(W_{2 1})\).

  The idea of the proof is to set up a diffeomorphism between \(W_{1
    2}\) and \(W_{2 1}\).  Our assumptions on \(C^x\)\hyp{}maps say
  that the resulting map on sections is a diffeomorphism.  Finally, we
  show that this diffeomorphism is the transition function defined in
  the statement of this proposition.

  Let \(\theta_1 \colon W_{1 2} \to T M\) be the map:
  \[
  \theta_1(t,v) = (\pi \times \eta_2)^{-1}(\alpha_2(t), \eta_1(v)).
  \]
  The definition of \(W_{1 2}\) ensures that \((\alpha_2(t),
  \eta_1(v)) \in V_2\) for \((t,v) \in W_{1 2}\) and this is the image
  of \(\pi \times \eta_2\).  Hence \(\theta_1\) is well-defined.  Define
  \(\theta_2 \colon W_{2 1} \to T M\) similarly.  These are both smooth maps.

  Notice that \(\pi(\pi \times \eta_i)^{-1} \colon V_i \subseteq U_i \times
  U_i \to U_i\) is the projection onto the first factor and \(\eta_i(\pi
  \times \eta_i)^{-1} \colon V_i \to U_i\) is the projection onto the second.
  Thus \(\pi \theta_1(t,v) = \alpha_2(t)\).  Hence \(\theta_1 \colon W_{1
    2} \to T M\) is such that \((t, \theta_1(t,v)) \in {\alpha_2}^* T
  M\) for all \((t,v) \in W_{1 2}\).  Then:
  \[
  \big(\alpha_1(t), \eta_2(\theta_1(t,v))\big) = (\alpha_1(t), \eta_1(v))
  \in V_1
  \]
  so \((t, \theta_1(t,v)) \in W_{2 1}\).  Hence we have a map
  \(\phi_{1 2} \colon W_{1 2} \to W_{2 1}\) given by:
  \[
  \phi_{1 2}(t,v) = (t, \theta_1(t,v)).
  \]
  Similarly we have a map \(\phi_{2 1} \colon W_{2 1} \to W_{1 2}\).  These
  are both smooth since the composition with the inclusion into \(S^1
  \times T M\) is smooth.

  Consider the composition \(\phi_{2 1}\phi_{1 2}(t,v)\).  Expanding
  this out yields:
  \begin{align*}
    \phi_{2 1}\phi_{1 2}(t,v) &= \phi_{2 1}(t, \theta_1(t,v)) \\
    &= (t, \theta_2(t,\theta_1(t,v))) \\
    &= \big(t, (\pi \times \eta_1)^{-1}(\alpha_1(t),
    \eta_2(\theta_1(t,v))) \big) \\
    &= \big(t, (\pi \times \eta_1)^{-1}(\alpha_1(t), \eta_1(v)) \big)
    \\
    &= \big(t, (\pi \times \eta)^{-1}(\pi(v), \eta_1(v)) \big) \\
    &= (t, v).
  \end{align*}
  The penultimate line is because \((t,v) \in {\alpha_1}^* T M\) so
  \(\pi(v) = \alpha_1(t)\).

  Hence \(\phi_{2 1}\) is the inverse of \(\phi_{1 2}\) and so
  \(\phi_{1 2} \colon W_{1 2} \to W_{2 1}\) is a diffeomorphism.  Thus by
  lemma~\ref{lem:cxsec}, the map \({\phi_{1 2}}_*\) is a
  diffeomorphism from \({\Psi_1}^{-1}(U_{1 2})\) to
  \({\Psi_2}^{-1}(U_{1 2})\).  We just need to show that this is
  the transition function.  It is sufficient to show that \(\Psi_2
  {\phi_{1 2}}_* = \Psi_2 \Phi_{1 2}\).  The right-hand side is, by
  definition, \(\Psi_1\), which satisfies:
  \[
  \Psi_1(\gamma)(t) = \eta_1(\tilde{\gamma}(t))
  \]
  where \(\tilde{\gamma} \colon S^1 \to T M\) is such that \(\gamma(t)
  = (t, \tilde{\gamma}(t))\).  On the other side:
  \begin{align*}
  {\phi_{1 2}}_*(\gamma)(t) &= \phi_{1 2}(\gamma(t)) \\
  &= \big(t, \theta_1(t, \tilde{\gamma}(t))\big) \\
  &= \big(t, (\pi \times \eta_2)^{-1}(\alpha_2(t),
  \eta_1(\tilde{\gamma}(t))) \big),
  \intertext{hence:}
  \Psi_2 {\phi_{1 2}}_*(\gamma)(t) &= \eta_2(\pi \times
  \eta_2)^{-1}( \alpha_2(t), \eta_1(\tilde{\gamma}(t))) \\
  &= \eta_1(\tilde{\gamma}(t)),
  \end{align*}
  as required.  Thus \({\phi_{1 2}}_* = \Phi_{1 2}\) and so the
  transition functions are diffeomorphisms.
\end{proof}

We therefore have a smooth atlas for \(L^x M\).

\section{Topology}
\label{sec:topology}

Following \cite[ch 27]{akpm} we proceed to topologise \(L^x M\) with
the inductive topology for the chart maps.  Our concern now is to
determine some topological properties of \(L^x M\).  The key is
theorem~\ref{th:submfd}.  Once we have proved this then the passage
from \(L^x \R\) to \(L^x M\) is straightforward.  We also show that
the inclusion \(L M \to L^x M\) is a homotopy equivalence.

\subsection{Submanifolds}

\begin{proposition}
Let \(M, N\) be finite dimensional smooth manifolds with \(M\) an
embedded submanifold of \(N\).  Then \(L^x M\) is an embedded
submanifold of \(L^x N\).
\end{proposition}

\begin{proof}
By the tubular neighbourhood theorem there is an open neighbourhood
\(U\) of \(M\) in \(N\), a smooth vector bundle \(E \to M\), and a
diffeomorphism \(\phi \colon E \to U\) which maps the zero section to the
inclusion of \(M\) in \(U\).  Let \(\eta \colon T M \to M\) be a local
addition over \(M\) with neighbourhood \(V \subseteq M \times M\).
Let \(E_V \subseteq E \times E\) be the restriction of \(E \times E\)
to \(V\).  Choose a connection on \(E\).  Let \((u,v) \in V\).  Let
\(p = (\pi \times \eta)^{-1}(u,v)\).  This lies in \(T_u M\) so the
path \(t \mapsto t p\) goes from the zero vector in \(T_u M\) to
\(p\).  Applying \(\eta\) results in a path from \(\eta(0_u) = u\) to
\(\eta(p) = v\).  Let \(P(u,v) \colon E_u \to E_v\) be the operator defined
by parallel transport along this path.

Using the connection a point in \(T E\) can be thought of as a
quadruple \((u,p,v,w)\) where \(u \in M\), \(p \in T_u M\), and \(v, w
\in E_u\) with the projection \(T E \to E\) being \((u,p,v,w) \to
(u,v)\) (we include \(u\) in the notation to emphasise the fibre).
Define \(\eta^E \colon T E \to E\) by \(\eta^E(u,p,v,w) =
(\eta(p),P(u,\eta(p))(v + w))\).  Then \(\pi^E \times \eta^E \colon T E \to
E \times E\) is:
\[
  (\pi^E \times \eta^E)(u,p,v,w) =
  \big((u,v),(\eta(p),P(u,\eta(p))(v + w)) \big).
\]
The projection of this to \(M \times M\) is \((u,\eta(p))\); so the
image of \(\pi^E \times \eta^E\) is in \(E_V\).  Since the map \((u,p)
\to (u,\eta(p))\) is onto \(V\), varying \(v\) and \(w\) shows that we
can get all of \(E_V\).  The inverse map is:
\[
  \big((u,v),(x,w)\big) \mapsto ((\pi \times \eta)^{-1}(u,x), v,
  P(\pi \times \eta^{-1}(u,x))(w) - v).
\]
This is smooth, so \(\pi^E \times \eta^E\) is a diffeomorphism onto
\(E_V\) and thus defines a local addition over the whole of \(E\). 

Using the diffeomorphism \(E \cong U\) we transfer this to \(U\) and
so get a local addition for \(M \subseteq N\).  The charts that this
defines make up part of the smooth atlas for \(N\).  Taking a chart
based at a smooth loop \(\alpha\) in \(M\), we see that the inclusion
of \(L^x M\) in \(L^x N\) looks like the inclusion of
\(\Gamma_{S^1}^x(\alpha^* T M)\) in \(\Gamma_{S^1}^x(\alpha^* T N)\).
This, in turn, looks like the inclusion of \(L^x \R^k\) in
\(L^x \R^n\).  Hence \(L^x M\) is an embedded submanifold
of \(L^x N\).
\end{proof}

\begin{corollary}
Let \(M\) be a closed finite dimensional smooth manifold.  Then the
following properties hold for \(L^x M\) if they hold for
\(L^x \R^n\): separable, metrisable, Lindel\"of, paracompact,
normal, smoothly regular, smoothly paracompact, and smoothly normal.
\end{corollary}

\begin{proof}
There is an embedding of \(M\) as a submanifold of \(\R^n\).  As it is
compact, the image is closed in \(\R^n\).  We therefore have an
embedding of \(L^x M\) in \(L^x \R^n\) which is also
closed.  The properties listed are all inherited by closed subspaces.
\end{proof}

\subsection{Homotopy Equivalence}
\label{sec:htpy}

One remarkable fact about the spaces \(L M\) and \(L^0 M\) is that
they are homotopy equivalent.  The standard method of this is to find
\emph{mollifiers} which ``smooth out'' continuous functions.  This
approach does not work with an arbitrary family of maps.  The paths
defined by the homotopy lie in the space of smooth loops at all points
except one end-point, therefore if smooth loops are not dense in the
given family of maps this homotopy cannot be continuous.  However
using the fact that \(L^x M\) is a smooth manifold one can still show
that the inclusion \(L M \to L^x M\) is a homotopy
equivalence.

The first step is to define the reverse map.  The basic idea is to use
a mollifier but we have to be a bit selective.  Let \(M\) be a closed
smooth finite dimensional manifold.  Via an embedding, regard \(M\) as
a submanifold of some Euclidean space, \(\R^n\).  By the tubular
neighbourhood theorem, there is a neighbourhood of \(M\) in \(\R^n\)
which retracts onto \(M\).  That is, there is some open neighbourhood
\(U \subseteq \R^n\) of \(M\) and a map \(p \colon U \to M\) which is
the identity on \(M\).  Let \(\eta \colon T M \to M \times M\) be a
local addition on \(M\) with image \(V\).

As \(M\) is compact we can find \(\mu > 0\) such that if \(x,y
\in M\) are such that \(\norm[x - y] < \mu\) then \((x,y) \in
V\).  We can also find \(\epsilon > 0\) such that if \(x \in M\) and
\(y \in \R^n\) are such that \(\norm[x - y] < \epsilon\) then \(y \in
U\) and \(\norm[x - p(y)] < \mu\).

\begin{lemma}
There is a continuous function \(\delta \colon L^0 \R^n \to \R\) such
that for \(\gamma \in L^0 \R^n\) then whenever \(\abs{s - t} <
\delta(\gamma)\), \(\norm[\gamma(s) - \gamma(t)] < \epsilon\).
\end{lemma}

\begin{proof}
Let \(\gamma \in L^0 \R^n\).  Then there is some \(\delta_\gamma
> 0\) such that whenever \(\abs{s - t} < \delta_\gamma\),
\(\norm[\gamma(s) - \gamma(t)] < \epsilon/3\).  Let \(\beta\) be such
that \(\norm[\beta - \gamma]_\infty < \epsilon/3\).  Then whenever
\(\abs{s - t} < \delta_\gamma\),
\begin{align*}
  \norm[\beta(s) - \beta(t)] &\le \norm[\beta(s) - \gamma(s)] +
  \norm[\gamma(s) - \gamma(t)] + \norm[\beta(t) - \gamma(t)]\\
&\le
  2 \norm[\beta - \gamma]_\infty + \norm[\gamma(s) - \gamma(t)] \\
&<
  \epsilon.
\end{align*}

Now \(L^0 \R^n\) is metrisable, hence paracompact and Hausdorff.
It therefore admits partitions of unity.  Choose a partition,
\(\{\tau_\lambda \colon \lambda \in \Lambda\}\), subordinate to the cover
of open balls of radius \(\epsilon/3\).  For each \(\lambda \in
\Lambda\) choose \(\gamma_\lambda\) such that the support of
\(\tau_\lambda\) is within the \(\epsilon/3\)\hyp{}ball around
\(\gamma_\lambda\).  Let \(\delta_\lambda \coloneqq
\delta_{\gamma_\lambda}\).  Define \(\delta \colon L^0 \R^n \to \R\)
by:
\[
  \delta(\gamma) = \sum_{\lambda \in \Lambda} \delta_\lambda
  \tau_\lambda(\gamma).
\]

For \(\gamma \in L^0 \R^n\), consider the set \(\Lambda(\gamma)
\coloneqq \{\lambda : \tau_\lambda(\gamma) \ne 0\}\).  This set is finite and
if \(\lambda \in \Lambda(\gamma)\) then \(\norm[\gamma -
  \gamma_\lambda] < \epsilon/3\).  As \(\delta(\gamma)\) is a convex
sum of the set \(\{\delta_\lambda : \lambda \in \Lambda(\gamma)\}\),
there is some \(\lambda \in \Lambda(\gamma)\) with \(\delta(\gamma)
\le \delta_\lambda\).  Whereupon we have that if \(\abs{s - t} <
\delta(\gamma)\), \(\norm[\gamma(s) - \gamma(t)] < \epsilon\) as
required.
\end{proof}

Using this, we define a continuous map \(R \colon L^0 \R^n \to
L \R^n\) with the property that \(\norm[\gamma -
  R(\gamma)]_\infty < \epsilon\) for all \(\gamma\).  Let \(\phi \colon
\R \to \R\) be a smooth bump function with support in \([-1,1]\) and
\(\int_\R \phi = 1\).  For \(r > 0\) let \(\phi_r \colon \R \to \R\) be the
function \(t \mapsto r \phi(t/r)\).  This has support in \([-r,r]\)
and satisfies \(\int_\R \phi_r = 1\).  Define \(R \colon L^0 \R^n \to
L \R^n\) by:
\[
  \gamma \mapsto \gamma * \phi_{\delta(\gamma)}.
\]
By this we mean that each component of \(\gamma\) is regarded as a map
with domain \R and is convoluted with the map
\(\phi_{\delta(\gamma)}\).

\begin{lemma}
The map \(R \colon L^0 \R^n \to L \R^n\) is continuous and
satisfies \(\norm[\gamma - R(\gamma)]_\infty < \epsilon\) for all
\(\gamma\).
\end{lemma}

\begin{proof}
Let \(\gamma \in L^0 \R\) and \(\beta \in \Ci(\R,\R)\).  Then:
\[
  (\gamma * \beta)(t) = \int_\R \gamma(s) \beta(t - s) d s.
\]
Hence:
\begin{align*}
  (\gamma * \beta)(t + 1) &= \int_\R \gamma(s) \beta(t + 1 - s) d s \\
&=
  \int_\R \gamma(\tilde{s} + 1) \beta(t - \tilde{s}) d \tilde{s} \\
&=
  \int_\R \gamma(\tilde{s}) \beta(t - \tilde{s}) d \tilde{s} \\
&= (\gamma
  * \beta)(t),
\end{align*}
whence \(R(\gamma)\) is periodic so can be viewed as a map with domain
\(S^1\).  The convolution of a continuous function by a smooth
function is again smooth so \(\gamma * \beta \in L \R\), with
derivative \(D (\gamma * \beta ) = \gamma * D\beta\).  Hence the image
of \(R\) is \(L \R^n\) as required.

To show that it is continuous, it is sufficient to show that the map
\(L^0 \R \times (0,\infty) \to L \R\), \((\gamma, r)
\mapsto \gamma * \phi_r\) is continuous.  Now for bounded maps on \R,
the map \((\alpha, \beta) \mapsto \alpha * \beta\) is bilinear and
satisfies:
\[
  \norm[\alpha * \beta]_\infty \le \norm[\alpha]_\infty
  \norm[\beta]_\infty
\]
where, by abuse of notation, we have used \(\norm_\infty\) for the
supremum norm for bounded functions on \R.  Hence for \(\alpha, \beta
\in L^0 \R\), \(r,s \in (0,\infty)\), and \(k \in \N\),
\[
  \norm[\alpha * (\phi_r)^{(k)} - \beta * (\phi_s)^{(k)}]_\infty \le
  \norm[\alpha]_\infty \norm[(\phi_r)^{(k)} - (\phi_s)^{(k)}]_\infty +
  \norm[\alpha - \beta]_\infty \norm[(\phi_s)^{(k)}]_\infty.
\]
This shows that providing \(\alpha\) and \(\beta\) are close and
providing \((\phi_r)^{(k)}\) and \((\phi_s)^{(k)}\) are close then
\((\alpha * \phi_r)^{(k)}\) and \((\beta * \phi_r)^{(k)}\) are close.
Thus to show that \(R\) is continuous it is sufficient to observe that
the map \(r \to \phi_r\) is continuous as a path into \(\Ci(\R,\R)\)
and this is straightforward.

Finally, let \(\gamma \in L^0 \R^n\).  For \(t \in S^1\), as
\(\int_\R \phi_r = 1\), \(\gamma(t) - R(\gamma)(t)\) is given by:
\[
  \gamma(t) - \int_\R \gamma(s) \phi_r(t - s) d s = \int_\R (\gamma(t)
  - \gamma(s)) \phi_{\delta(\gamma)}(t - s) d s.
\]
Now \(\phi_{\delta(\gamma)}(t - s)\) is zero outside \([t -
  \delta(\gamma), t + \delta(\gamma)]\) and on this interval
\(\abs{\gamma(t) - \gamma(s)} < \epsilon\), by definition of
\(\delta(\gamma)\).  Hence \(\abs{\gamma(t) - R(\gamma)(t)} <
\epsilon\) as required.
\end{proof}

\begin{corollary}
\label{cor:manmoll}
There is a continuous map \(R_M \colon L^0 M \to L M\) with
the property that for all \(\gamma \in L^0 M\), \((\gamma(t),
R_M(\gamma)(t)) \in V\) for all \(t \in S^1\).
\end{corollary}

\begin{proof}
We restrict the map \(R\) to the domain \(L^0 M\).  By
construction, \(R(\gamma)\) takes values in \(U\), the neighbourhood
of \(M\).  The map \(R_M\) is the composition of this with the
projection \(p \colon U \to M\).  The required property holds because of
the choices made.
\end{proof}

\begin{theorem}
Let \(L^x \R\) be a class of maps satisfying the assumptions of
section~\ref{sec:intro}.  Then the inclusion \(L M \to
L^x M\) is a homotopy equivalence.
\end{theorem}

\begin{proof}
The reverse map is the composition of the inclusion of \(L^x M\)
in \(L^0 M\) with the map \(R_M\).  We will denote this again by
\(R_M\).

By construction, for \(\gamma \in L^x M\), \((R_M(\gamma),
\gamma)\) takes values in \(V\).  Define \(H \colon L^x M \times
[0,1] \to L^x M\) by:
\[
  H(\gamma,s) = \pi_2 \eta(s \eta^{-1}(R_M(\gamma), \gamma)),
\]
where in \(T M\) we have used the natural \(\R\)\hyp{}action on the
fibres and \(\pi_2\) is the projection onto the second factor.  This
is continuous as it is the composition of continuous maps.  For \(s =
1\) we have \(H(\gamma,1) = \pi_2(R_M(\gamma),\gamma) = \gamma\).  For
\(s = 0\), \(H(\gamma,0) = \pi_2(\eta(R_M(\gamma), R_M(\gamma))) =
R_M(\gamma)\).  Hence \(H\) is the required homotopy.

The same homotopy works for smooth maps, showing that the composition
\(L M \to L^x M \to L M\) is homotopic to the identity.
\end{proof}

\subsection{Based Loops}
\label{sec:based}

All of our discussion so far holds for based loops as well as free
loops.  For the main part, the only modification needed is the
insertion of the word ``based'' at appropriate points.  The only place
where more is required is in proving the homotopy
equivalence.  The problem there is that the result of applying a
mollifier to a based continuous loop may no longer be based.  The fix
is simple, however, since at that point in the construction of the
homotopy equivalence we are dealing with loops in \(\R^n\).  We can
therefore define the based mollifier \(R_0\) in terms of the map \(R\)
defined in section~\ref{sec:htpy} as
\[
  R_0(\gamma) = R(\gamma) - R(\gamma)(0)
\]
(we are tacitly assuming that the basepoint of \(\R^n\) is the
origin).  To ensure that the resulting map \(R_0(\gamma)\) has the
properties analogous to corollary~\ref{cor:manmoll} we need to replace
\(\epsilon\) by \(\epsilon/2\) in section~\ref{sec:htpy}.

The relationship between based loops and free loops is an important
one.  In homotopy theory there is a fibration
\[
  \Omega^0 M \to L^0 M \to M
\]
and so, via the homotopy equivalences, we can deduce that
\[
  \Omega^x M \to L^x M \to M
\]
is also a fibration.  This describtion, however, is not one of the
flavour of differential topology.  A more suitable description would
be that it is a locally trivial fibre bundle.  This will follow from
the following theorem.  Let \(e_0 \colon L^x M \to M\) be the
evalutation map at \(0\), \(\gamma \mapsto \gamma(0)\).

\begin{theorem}
Let \(M, P\) be a smooth finite dimensional orientable manifolds
without boundary.  Suppose that \(P\) is an embedded submanifold of
\(M\) with tubular neighbourhood \(U \subseteq M\) and normal bundle
\(N \to P\).  Define
\[
  L^x_P M \coloneqq \{\gamma \in L^x M : \gamma(0) \in P\}
\]
and \(L^x_U M\) similarly.  Then \(L^x_P M\) is an embedded
submanifold of \(L^x M\) with tubular neighbourhood \(L^x_U M\) and
normal bundle \({e_0}^* N\).  Moreover, the evalutation map \(e_0
\colon L^x M \to M\) takes the quadruple \((L^x M, L^x_P M ,L^x_U M,
{e_0}^* N)\) to \((M, P, U, N)\) preserving all the structure.
\end{theorem}

\begin{proof}
We omit the proof that \(L^x_P M\) is a smooth submanifold of \(L^x
M\) as this is a simple modification of the work of previous sections;
the model space is \(L^x_{\R^k} \R^n\).  The case of \(L^x_U M\) is
simpler as it is an open submanifold of \(L^x M\).

Let \(\pi \colon N \to P\) be the projection and let \(\phi \colon U
\to N\) be the diffeomorphism.  Our strategy is to find a continuous
map \(\Psi \colon U \to \diff_c(U)\), where \(\diff_c(U)\) is
diffeomorphisms of \(U\) that are the identity outside a compact set,
such that \(\Psi_u(u) = \pi\phi(u)\).  Using this, we define the
diffeomorphism \(L_U^x M \to {e_0}^* N\) by
\[
  \alpha \mapsto \left(\Psi_{\alpha(0)}(\alpha), \phi(\alpha(0))\right)
\]
with inverse
\[
  (\beta,v) \mapsto \Psi_{\phi^{-1}(v)}^{-1}(\beta).
\]
That these are smooth follows from the fact that they are defined
entirely in terms of smooth maps of the original manifolds and these
induce smooth maps on our loop spaces by assumption.

Thus we need to find the map \(\Psi \colon U \to \diff_c(U)\) with the
appropriate conditions.  We actually define the map for \(N\) as there
we can use the vector bundle structure and then transfer it to \(U\)
via the diffeomorphism \(\phi\).

The first stage of defining \(\Psi\) is to define a map \(s \colon N
\to \Gamma_c(N)\), the space of sections of \(N\) with compact
support.  Let \(\{(U_\lambda, \nu_\lambda, \rho_\lambda) : \lambda \in
  \Lambda\}\) be a family of triples where
\begin{enumerate}
\item \(\{U_\lambda : \lambda \in \Lambda\}\) is a locally finite open
cover of \(P\) such that each \(U_\lambda\) has compact closure in
\(P\);
\item \(\nu_\lambda \colon \pi^{-1}(U_\lambda) \to U_\lambda \times
\R^n\) is a trivialisation of \(N\) over \(U_\lambda\);

\item \(\{\rho_\lambda : \lambda \in \Lambda\}\) squares to a
partition of unity subordinate to the cover \(\{U_\lambda\}\); that
is, \(\rho_\lambda \colon P \to \R\) is a bump function with support
in \(U_\lambda\) and \(\sum_\lambda \rho_\lambda(x)^2 = 1\) for
all \(x \in P\).
\end{enumerate}

Let \(\tilde{\nu}_\lambda \colon \pi^{-1}(U_\lambda) \to \R^n\) be the
composition of \(\nu_\lambda\) with the projection onto \(\R^n\).
Define \(s \colon N \to \Gamma_c(N)\) by
\[
  s(v)(x) \coloneqq \sum_\lambda \rho_\lambda(\pi(v)) \rho_\lambda(x)
  \nu_\lambda^{-1}(x,\tilde{\nu}_\lambda(v)).
\]
The sum is well-defined since \(\rho_\lambda(\pi(v))\rho_\lambda(x)\)
can only be non-zero if both \(\pi(v)\) and \(x\) are in the domain of
\(\nu_\lambda\).  Each \(s(v)\) is clearly a section and its support
is contained in the finite union \(\bigcup \{U_\lambda : \pi(v) \in
  U_\lambda\}\), and hence has compact support.  Observe that
\[
  s(v)(\pi(v)) = \sum_\lambda \rho_\lambda(\pi(v))\rho_\lambda(\pi(v))
  \nu_\lambda^{-1}(\pi(v),\tilde{\nu}_\lambda(v)) = \sum_\lambda
  \rho_\lambda(\pi(v))^2 v = v.
\]

There is a canonical embedding of \(N\) in \(T N\) as the vertical
tangent bundle and thus we can extend any section \(\sigma\) of \(N\)
to a vector field on \(N\) by defining \(\tilde{\sigma}(v) =
\sigma(\pi(v))\).  If the original section had compact support then
the resulting vector field will have compact horizontal support.  We
wish to apply this proceedure to the sections that we have defined
above, but we also wish to ensure that the resulting vector fields
have genuine compact support.  To do that, we choose an inner product
on the fibres of \(N\) which varies smoothly over \(P\) and a bump
function \(\tau \colon \R \to \R\) which takes the value \(1\) on \([
  0, 1 ]\) and is zero above, say, \(2\).  Define
\[
  X_v(u) \coloneqq - \tau \left(\norm[u]^2/(1 + \norm[v]^2)\right)
  s(v)(\pi(u)).
\]
This has compact support, is a vertical vector field, and for \(u \in
N_v\) with \(\norm[u] \le \norm[v]\) then \(X_v(u) = - v\).

We therefore have a continuous map \(N \to \Xi_c(N)\).  We compose
this with the exponential map \(\exp \colon \Xi_c(N) \to \diff_c(N)\)
to define \(\Psi \colon N \to \diff_c(N)\).  The properties of \(X_v\)
translate into properties of \(\Psi_v\).  As \(X_v\) is a vertical
vector field, \(\Psi_v\) preserves the fibres of \(N\).  Most
importantly, \(\Psi_v(v) = 0_v\).

This is the required map and so establishes \(L^x_U M\) as the tubular
neighbourhood of \(L^x_P M\).  It is clear from the setup that the
evaluation map has the properties stated in the theorem.
\end{proof}

\begin{corollary}
The evaluation map \(L^x M \to M\) is a locally trivial fibre bundle
with fibre \(\Omega^x M\).
\end{corollary}

\begin{proof}
Take \(P = \{x_0\}\) to be the basepoint and \(U\) the codomain of a
chart near \(x_0\) with domain \(\R^n\).
\end{proof}

\section{Circle Actions}
\label{sec:diff}

The diffeomorphism group of the circle acts on maps with domain
\(S^1\) by precomposition.  It is usual to assume that \(L^x \R\) is
closed under this action, whence we get an action on \(L^x M\) for
\(M\) a smooth finite dimensional manifold.  We would like to transfer
knowledge of that action from \(L^x \R\) to \(L^x M\).

\subsection{Transferring the Action}

We are going to prove an inheritance result which states that the
action on \(L^x M\) is the same as that on \(L^x \R^n\).  This will be
an easy corollary of theorem~\ref{th:submfd}.  The more important
point of this section is to consider what types of action there are.

\begin{defn}
Let \(M\) be a smooth finite dimensional manifold.  Let \(G \subseteq
\diff(S^1)\) be a sub\hyp{}Lie group of the set of diffeomorphisms of
the circle.  We define the following possible types of action of \(G\)
on \(L^x M\).

\begin{enumerate}
\item The action is by bijections.

\item The action is by homeomorphisms.

\item The action is by diffeomorphisms.

\item The action map, \(G \times L^x M \to L^x M\), is continuous.

\item The action map, \(G \times L^x M \to L^x M\), is smooth.

\item The representation map, \(G \to \homeo(L^x M)\), is
continuous.

\item The representation map, \(G \to \diff(L^x M)\), is
smooth.
\end{enumerate}
\end{defn}

The diffeomorphism group, \(\diff(S^1)\), is an open subset of \(L
S^1\) and thus inherits the structure of a smooth manifold.  The
circle, acting by rigid rotations, is a subset of \(\diff(S^1)\) and
the inherited structure is the same as its usual one.

Although these levels have been written in a necessarily linear form,
the relationships between them are more complicated than this
suggests.  For example, a continuous representation map does not
necessarily imply a continuous action map as the evaluation map
\(\homeo(X) \times X \to X\) is not necessarily jointly continuous.

\begin{proposition}
All the levels defined are inherited by \(L^x M\) from \(L^x \R\).
\end{proposition}

\begin{proof}
As diffeomorphisms of the circle act linearly on \(L^x \R\) this
proposition holds for \(M = \R^n\) simply by taking finite products.
The result for general \(M\) follows from theorem~\ref{th:submfd}.
Since \(L^x M\) is an embedded submanifold of \(L^x \R^n\), a map into
\(L^x M\) is continuous or smooth if and only if it is continuous or
smooth into \(L^x \R^n\); and the restriction to \(L^x M\) of a
continuous or smooth map from \(L^x \R^n\) is again continuous or
smooth.

For the representation maps, we are being deliberately vague about the
topologies on the homeomorphism and diffeomorphism groups.  There is
considerable freedom in choosing this topology and we wish to allow
for this freedom, only assuming that the topologies are compatbile for
\(M\) as for \(\R^n\).  Then as \(L^x M\) is a smooth retract of an
open subset of \(L^x \R^n\), the restriction map from the subspace of
homeomorphisms, resp.\ diffeomorphisms, of \(L^x \R^n\) which preserve
\(L^x M\) as a set to the set of homeomorphisms, resp.\
diffeomorphisms, of \(L^x M\) is defined and continuous, resp.\
smooth.  As the representation map of \(G\) factors through this map
its properties are inherited from those of the representation on \(L^x
\R\).
\end{proof}

\subsection{Continuity of Linear Circle Actions}

The most common action considered on loop spaces is that of the circle
itself.  In light of the inheritance properties of circle actions, it
seems a good idea to consider the general case of the circle acting on
a locally convex topological vector space.  As this is, by its very
nature, more in the realm of functional analysis than differential
topology, at each stage we shall consider how it applies to the
examples of smooth loops and continuous loops in order to ground the
discussion in terms familiar to the differential topologist.

We start with a more detailed discussion of what it might mean for a
circle action to be ``continuous''.  There are several ``levels'' of
continuity that one could consider, more than those listed in the
previous section.  The following definition contains the ones that we
think are interesting or useful.

\begin{defn}
\label{def:cts}
Let \(E\) be a lctvs.  Suppose that the circle acts on \(E\) by linear
maps, not necessarily continuous.  Let \(R_t\) be the linear map
corresponding to \(t \in S^1\).  We define the following levels of
continuity for this action:
\begin{enumerate}
\item The representation is continuous; that is, the action induces a
continuous map \(S^1 \to \m{L}_b(E)\).  Here, \(\m{L}_b(E)\) denotes
the space of continuous linear maps from \(E\) to itself equipped with
the \emph{strong} topology; that is, the topology of uniform
convergence on bounded sets.
\label{it:rep}

\item The action is continuous; that is, the action is continuous as a
map \(S^1 \times E \to E\).
\label{it:cts}

\item The action is separately continuous; that is, for each \(t \in
S^1\) then \(x \to R_t x\) is a continuous map \(E \to E\), and
for each \(x \in E\) then \(t \to R_t x\) is a continuous map \(S^1
\to E\).
\label{it:sep}

\item The action is by equicontinuous linear maps; that is to say, for
each \(0\)\hyp{}neighbourhood \(V\) in \(E\) there is a
\(0\)\hyp{}neighbourhood \(U\) such that \(R_t U \subseteq V\) for all \(t
\in S^1\).
\label{it:equi}

\item There is a \(0\)\hyp{}basis of \(S^1\)\hyp{}invariant sets.
\label{it:nbdin}

\item The action is by continuous linear maps; that is, each \(R_t\)
is continuous.
\label{it:clin}

\item The topology on \(E\) is \(S^1\)\hyp{}invariant.
\label{it:topin}
\end{enumerate}
\end{defn}

The strong topology is the finest topology that one would sanely use.
Thus positive results for the strong topology will propagate forwards
to any coarser topology.  A \(0\)\hyp{}basis for this topology
consists of the sets:
\[
  N(B,U) \coloneqq \{T \in \m{L}(E) : T(B) \subseteq U\}
\]
where \(B, U\) are subsets of \(E\) with \(B\) bounded and \(U\) a
\(0\)\hyp{}neighbourhood.  If \(E\) is a Banach space then this is the
usual topology which is normable with norm \(\norm[T] \coloneqq \sup\{\norm[T
    x] : \norm[x] \le 1\}\)..

We shall now show how the list in definition~\ref{def:cts} is,
roughly, from the strictest to the weakest.  We start with just the
results that apply to all lctvs.

\begin{proposition}
Let \(E\) be a lctvs with an action of the circle by linear maps.  We
have the following links between the levels of continuity:
\begin{enumerate}
\renewcommand{\theenumi}{(\roman{enumi})}
\item \ref{it:equi} is equivalent to \ref{it:nbdin};
\label{it:eqnbd}

\item \ref{it:clin} is equivalent to \ref{it:topin};
\label{it:ctop}

\item \ref{it:cts} is equivalent to having both \ref{it:sep} and
\ref{it:equi};
\label{it:cteq}

\item \ref{it:sep} implies \ref{it:clin};
\label{it:secl}

\item \ref{it:rep} implies \ref{it:sep}.
\label{it:reep}
\end{enumerate}
\end{proposition}

Before proving this we remark that the reason why \ref{it:rep} does
not automatically imply \ref{it:cts} is because the evaluation map
\(\m{L}_b(E) \times E \to E\) is not, in general, continuous but only
separately continuous.  Thus the action map is separately continuous
as it factors as:
\[
  S^1 \times E \to \m{L}_b(E) \times E \to E
\]
but we cannot deduce from this that it is continuous.

\begin{proof}
The equivalences \ref{it:eqnbd} and \ref{it:ctop} are obvious, as
is the implication \ref{it:secl}.  The deduction of \ref{it:sep}
from \ref{it:cts} is also obvious.  We have already explained
\ref{it:reep}.

Thus only \ref{it:cteq} remains and of that we need to show that
\ref{it:cts} implies \ref{it:equi} and that together
\ref{it:sep} and \ref{it:equi} imply \ref{it:cts}.

To show that \ref{it:cts} implies \ref{it:equi} let \(V\) be an
open \(0\)\hyp{}neighbourhood in \(E\).  By assumption, for each \(t \in
S^1\) there is some open \(0\)\hyp{}neighbourhood \(U_t\) and \(\delta_t >
0\) such that \((t - \delta_t, t + \delta_t) \times U_t\) maps into
\(V\).  As \(S^1\) is compact there is some finite set \(\{t_1,
  \dotsc, t_n\}\) such that the intervals \(\{(t_j - \delta_{t_j}, t_j
  + \delta_{t_j})\}\) cover \(S^1\).  Let \(U = \bigcap_{j=1}^n
U_{t_j}\).  Then \(U\) is a finite intersection of open
\(0\)\hyp{}neighbourhoods, hence is one itself.  For \(t \in S^1\) there
is some \(j\) such that \(t \in (t_j - \delta_{t_j}, t_j +
\delta_{t_j})\) whence, as \(U \subseteq U_{t_j}\), \(R_t(U) \subseteq
V\).  Thus the action is by equicontinuous linear maps.

For the converse we assume both \ref{it:sep} and \ref{it:equi}.
Let \(x \in E\) and \(t \in S^1\).  Let \(V\) be a convex
\(0\)\hyp{}neighbourhood which, by \ref{it:equi}, we may assume to be
\(S^1\)\hyp{}invariant.  Then \(\frac12V\) is also a convex,
\(S^1\)\hyp{}invariant \(0\)\hyp{}neighbourhood so as the map \(s \to R_s x\)
is continuous at \(t\) there is some \(\delta > 0\) such that if
\(\abs{s} < \delta\) then \(R_t x - R_{t+s} x \in \frac12 V\).  Let
 \(s \in S^1\) be such that \(\abs{s} < \delta\) and let \(y \in x +
\frac12 V\).  We have:
\begin{align*}
R_t x - R_{t+s} y &= R_t x - R_{t+s} x + R_{t+s} x - R_{t+s} y \\
&= R_t x - R_{t+s} x + R_{t+s}(x - y) \\
&= \frac12 (2 R_t x - 2 R_{t+s} x) + \frac12 R_{t+s}(2 x - 2 y).
\end{align*}
Now \(2 R_t x - 2 R_{t+s} x\) and \(2 x - 2 y\) both lie in \(V\).  As
\(V\) is \(S^1\)\hyp{}invariant, \(R_{t+s}(2 x - 2 y)\) also lies in
\(V\).  Thus as \(V\) is convex, \(R_t x - R_{t+s} y\) is in \(V\).
Hence \((t - \delta, t+ \delta) \times x + \frac12 V\) lies in the
preimage of \(R_t x + V\).  Hence the action is continuous.
\end{proof}

There are more connections between these conditions if the space \(E\)
has more structure.

As mentioned above the failure of \ref{it:rep} to automatically
imply \ref{it:cts} is due to possibility that the evaluation map is
not continuous.  It is continuous if, and only if, \(E\) is normable.
Thus we deduce:

\begin{lemma}
\label{lem:normcts}
Let \(E\) be a normable lctvs with an action of the circle by linear
maps.  Then \ref{it:rep} implies \ref{it:cts}. \hspace*{\fill}\qedsymbol
\end{lemma}

A more general class of spaces that allows us to strengthen the links
is the family of \emph{barrelled} lctvs.  This is a technical property
of lctvs which we shall not describe here, we merely need one of its
well-known consequences.  It follows from \cite[11.1.5]{hj} and
Baire's theorem that \(L^0 \R\) and \(L \R\) are barrelled.

\begin{proposition}
\label{prop:barrel}
Let \(E\) be a barrelled lctvs with an action of the circle by linear
maps.  Then \ref{it:sep} implies \ref{it:equi}.  Hence each of
\ref{it:rep} and \ref{it:sep} imply \ref{it:cts}.
\end{proposition}

\begin{proof}
The proof that \ref{it:sep} implies \ref{it:equi} is similar to
\cite[III\S5.3]{hs}.  That the action is separately continuous means
that the map \(S^1 \to \m{L}(E)\) is well-defined and is continuous
for the topology of uniform convergence on all \emph{finite} sets.
Thus the image of \(S^1\) in \(\m{L}(E)\) is simply bounded and hence,
as \(E\) is barrelled, equicontinuous.

Since \ref{it:sep} and \ref{it:equi} together imply \ref{it:cts}
we therefore have that \ref{it:sep} alone implies \ref{it:cts}.
Also as \ref{it:rep} implies \ref{it:sep} we also have that
\ref{it:rep} implies \ref{it:cts}.
\end{proof}

A useful property of \(L \R\) is that closed bounded subsets are
compact; this follows from \cite[II\S7.2]{hs} as it is a complete
nuclear space.

\begin{proposition}
\label{prop:ctsrep}
Let \(E\) be a lctvs with an action of the circle by linear maps.
Suppose that every closed, bounded subset of \(E\) is compact.  Then
\ref{it:cts} implies \ref{it:rep}.
\end{proposition}

\begin{proof}
We shall show that if the action is continuous then the map \(S^1 \to
\m{L}(E)\) is continuous for the topology of uniform convergence on
compact sets.  The assumption on \(E\) then says that this is
precisely the topology of uniform convergence on bounded sets.

So assume that the circle action on \(E\) is continuous.  Let \(C, V
\subseteq E\) be such that \(C\) is compact and \(V\) is a convex,
circled \(0\)\hyp{}neighbourhood.  Let \(t_0 \in S^1\).  As the circle
action is continuous then for each \(c \in C\) there is some
\(\delta_c > 0\) and \(U_c\) a neighbourhood of \(c\) in \(E\) such
that if \(x \in U_c\) and \(\abs{t} < \delta_c\) then \(R_{t_0 + t} x
- R_t c \in \frac12V\).

The neighbourhoods \(\{U_c\}\) cover \(C\) so there is some finite
subset which will do; say, \(U_1, \dotsc, U_n\) corresponding to
points \(c_1, \dotsc, c_n \in C\).  Let \(\delta\) be the minimum of the
corresponding subfamily of \(\{\delta_c\}\);  then \(\delta > 0\).

Let \(t\) be such that \(\abs{t} < \delta\).  Let \(c \in C\), then
there is some \(j\) such that \(c \in U_j\).  Thus \(R_{t_O + t} c -
R_{t_0} c_j \in \frac12 V\).  Now the choice of \(c_j\) depended only
on \(c\) and not on \(t\).  Therefore we also have \(R_{t_O + 0} c -
R_{t_0} c_j \in \frac12 V\).  Thus:
\[
  R_{t_0 + t} c - R_{t_0} c = R_{t_0 + t} c - R_{t_0} c_j + R_{t_0}
  c_j - R_{t_0} c
\]
which, for the usual convexity reasons, lies in \(V\).  Hence for
\(\abs{t} < \delta\), \(R_{t_0 + t} - R_{t_0}\) maps \(C\) into
\(V\).  Thus the map \(S^1 \to \m{L}(E)\) is continuous for the
topology of uniform convergence on compact subsets.
\end{proof}

\subsection{Circle Actions on Loop Spaces}

In this section we shall use the results of the previous one to
determine how continuous are the circle actions on our example
spaces.  For convenience we list the technical properties of our
spaces so that we know which of the above results apply.  We also, for
quick reference, list a \(0\)\hyp{}basis.

\begin{enumerate}
\item \(L^0 \R\) is barrelled.  The topology is determined by
the sets:
\[
  U(\epsilon) \coloneqq \{\gamma : \sup\{\abs{\gamma(t)} : t \in S^1\} <
    \epsilon\}.
\]

\item \(L \R\) is barrelled and every closed bounded subset is
compact.  The topology is determined by the sets:
\[
  U(n,\epsilon) \coloneqq \{\gamma : \sup\{\abs{\gamma^{(k)}(t)} : t \in S^1,
      0 \le k \le n\}< \epsilon \}
\] 
\end{enumerate}

We shall now determine how continuous is the action of rotation of
loops on each of these spaces. 

\begin{proposition}
For both spaces the action is by continuous linear maps.
\end{proposition}

\begin{proof}
We just need to show that the topology is \(S^1\)\hyp{}invariant.  It is
sufficient to show this for the \(0\)\hyp{}neighbourhoods listed above.
We have:
\begin{gather*}
  R_s U(\epsilon) = U(\epsilon), \\
  R_s U(n, \epsilon) = U(n, \epsilon),
\end{gather*}
Thus in each case the topology is \(S^1\)\hyp{}invariant and so the action
is by continuous linear maps.
\end{proof}

\begin{proposition}
For both spaces the action is by equicontinuous linear maps.
\end{proposition}

\begin{proof}
From the previous proof it is obvious that the given \(0\)\hyp{}basis is of
\(S^1\)\hyp{}invariant sets.  Hence for these three the action is by
equicontinuous linear maps.
\end{proof}

\begin{proposition}
The circle action on both of \(L^0 \R\) and \(L \R\) is separately
continuous.
\end{proposition}

\begin{proof}
We already know that the circle acts by continuous linear maps which
is half of separate continuity.  Thus we need to show that for each
loop \(\gamma\) then the map \(t \to R_t \gamma\) is continuous.

We shall give the proof in full for \(L \R\).  The proof for
\(L^0 \R\) is a simplification of this.  We need to show that for
\(\gamma \in L \R\), \(t_0 \in S^1\), and a
\(0\)\hyp{}neighbourhood \(V\) then there is some \(\delta > 0\) such that
if \(\abs{t} < \delta\) then \(R_{t_0 + t} \gamma - R_{t_0} \gamma \in
V\).  It is sufficient to do this for \(V = U(n, \epsilon)\) whence
we need to show that \(\norm[(R_{t_0 + t} \gamma)^{(k)} - (R_{t_0}
  \gamma)^{(k)}]_\infty < \epsilon\) for \(0 \le k \le n\).

Expanding out the definition of the norm and using 
\((R_s \alpha)^{(k)}(t) = \alpha^{(k)}(t + s)\) we see that we want to
ensure that:
\[
  \sup\{\abs{\gamma^{(k)}(s + t) - \gamma^{(k)}(s)} : s \in S^1, 0 \le
    k \le n\} < \epsilon
\]
whenever \(\abs{t} < \delta\).  That such a \(\delta > 0\) exists
comes from the fact that the loops \(\gamma, \gamma^{(1)}, \dotsc,
\gamma^{(n)}\) are all uniformly continuous on \(S^1\) and there is
only a finite number of them.

For \(L^0 \R\) the situation is slightly simplified in that we
only need to consider \(\gamma\) and not any of its derivatives (which
it may not have, of course).
\end{proof}

From propositions~\ref{prop:barrel} and~\ref{prop:ctsrep} we deduce
the following.

\begin{corollary}
The circle actions on \(L^0 \R\) and \(L \R\) are continuous.  The
representation for the action on \(L \R\) is also continuous.  \hspace*{\fill}\qedsymbol
\end{corollary}

Thus the action on \(L \R\) is the best it can be.  This is not true
of \(L^0 \R\).  We deduce this from a more general result which says
that this is not the fault of the type of loop but rather of using a
normed vector space of loops.  Recall that a trigonometric polynomial
is a (finite) linear span of sines and cosines.

\begin{proposition}
Let \(E \subseteq \map(S^1, \R)\) be an \(S^1\)\hyp{}invariant vector
space of loops which contains the trigonometric polynomials.  Let
\(p\) be an \(S^1\)\hyp{}invariant semi-norm on \(E\) which restricts to a
norm on the subspace of trigonometric polynomials.  Let \((\tilde{E},
\norm)\) be the associated Banach space.  Then circle action is by
equicontinuous linear maps but the associated representation is not
continuous.
\end{proposition}

Thus in this general case the only question to answer is whether or not
the circle action itself is continuous.

\begin{proof}
As the set-up is \(S^1\)\hyp{}invariant, the unit ball in \(\tilde{E}\) is
\(S^1\)\hyp{}invariant and so the circle acts by equicontinuous linear
maps.

Let \(\delta > 0\).  Choose \(n \in \N\) such that \(1/n < \delta\).
As \(E\) contains the trigonometric polynomials it contains the loop
\(\gamma(t) = \cos(2 \pi n t) v\) where \(v \in \R\) is non-zero.
By assumption on the semi-norm, \(\gamma\) represents a non-zero
element in \(\tilde{E}\).  Let \(h = 1/(2n)\), then \(R_h \gamma = -
\gamma\).  Hence \(\norm[(I - R_h)\gamma] = 2 \norm[\gamma]\) and so
\(\norm[I - R_h] \ge 2\).  Thus the map \(t \to R_t\) is not
continuous into \(\m{L}_b(\tilde{E})\).
\end{proof}

In summary, the circle action on \(L \R\) is as good as it can be
whereas that on \(L^0 \R\) is almost that good and is as good as it
can be given that it is a normed vector spaces.

\subsection{Smooth Actions}

We conclude with a comment on how smooth are the circle actions on \(L
\R\) and on \(L^0 \R\).  As with continuity we can ask for different
levels of smoothness.

For the positive results in this section we have to decide on a type
of calculus.  We choose the convenient calculus of \cite{akpm}.  This
states that a map into a locally complete lctvs is smooth if and only
if its composition with each continuous linear functional is a smooth
map into \R.  This provides us with test functions to determine
whether or not a map is smooth.  For the negative results we do not
need to pick a calculus as for any calculus, continuous linear maps
are certainly smooth and so we can still use them as test functions to
determine if a map is not smooth.

We start with some positive results about \(L \R\).

\begin{proposition}
\label{prop:smoothissmooth}
The action map \(\rho \colon S^1 \times L \R \to L \R\) is
smooth.
\end{proposition}

\begin{proof}
As \(L \R\) is a closed subspace of \(\Ci(\R, \R)\) and
\R is a covering space of \(S^1\), it is clearly sufficient to show
that the map \(\tilde{\rho} \colon \R \times \Ci(\R, \R) \to \Ci(\R, \R)\),
\((t, \zeta) \to (s \mapsto \zeta(s + t))\), is smooth.  We need to
show that it takes smooth curves in \(\R \times \Ci(\R, \R)\) to
smooth curves in \(\Ci(\R, \R)\).

Let \(c \colon \R \to \R \times \Ci(\R, \R)\) be smooth.  Let \(\tilde{c} =
\tilde{\rho} \circ c\).  We can write \(c = (c_1, c_2)\) for smooth
curves \(c_1 \colon \R \to \R\) and \(c_2 \colon \R \to \Ci(\R, \R)\) since the
obvious projection maps are smooth.  Then for \(t \in \R\)
\[
  \tilde{c}(t) = (s \mapsto c_2(t)(s + c_1(t))).
\]
By the exponential law, \(\tilde{c}\) is smooth if and only if its
adjoint, \(\tilde{c}^\lor \colon \R^2 \to \R\) is smooth.  This adjoint is
\[
  (s,t) \mapsto c_2(t)(s + c_1(t)).
\]
Now as \(c_2 \colon \R \to \Ci(\R, \R)\) is smooth its adjoint,
\(c_2^\lor\), is also smooth, again by the exponential law.  This is
the map \((s,t) \mapsto c_2(t)(s)\).  Thus \(\tilde{c}^\lor\) is
smooth as it factors as the composition
\[
  (s,t) \mapsto (s + c_1(t),t) \xrightarrow{c_2^\lor} c_2(t, s +
  c_1(t)).
\]
Hence \(\tilde{\rho}\) maps smooth curves to smooth curves and is thus
smooth.
\end{proof}

\begin{corollary}
The representation map \(S^1 \to \m{L}_b(L \R)\) is smooth.
\end{corollary}

\begin{proof}
This follows from the uniform boundedness principle, see
\cite[I.5.18]{akpm}: a map into \(\m{L}_b(L \R)\) is smooth if and
only if all composites with evaluations at points in \(L \R\) are
smooth.
\end{proof}

Note that we cannot deduce from this that \(S^1 \to \m{L}_b(L
\R)\) is continuous since we have left the realm where the
\(\ci\)\hyp{}topology agrees with the locally convex topology one.

Now we turn to the negative result and recall that here we do not
assume a particular calculus.

\begin{proposition}
\label{prop:mustbesmooth}
Let \(L^x \R\) be a class of loops satisfying the conditions
\ref{cond:vspace}, \ref{cond:lctvs}, and \ref{cond:smthcts} of the
introduction.  Let \(\gamma \in L^x \R\) be such that the map \(S^1
\to L^x \R\), \(t \mapsto R_t \gamma\), is smooth.  Then \(\gamma
\colon S^1 \to \R\) is smooth.

If, in addition, for all \(\gamma \in L \R\) then the maps \(S^1 \to L
\R\), \(t \mapsto R_t \gamma\), are smooth then the above becomes an
if\hyp{}and\hyp{}only\hyp{}if.
\end{proposition}

\begin{proof}
Let \(e_0 \colon L^x \R \to \R\) be the evaluation map at \(0\).  This
is continuous by the assumptions and hence is smooth.  As \(S^1 \to
L^x \R\), \(t \mapsto R_t \gamma\), is smooth its composition with
\(e_0\) is a smooth map \(S^1 \to \R\).  This composition is \(t
\mapsto e_0(R_t \gamma) = \gamma(0 + t) = \gamma(t)\).  Thus
\(\gamma\) is smooth.

For the second part, let \(\gamma \in L^x \R\) be a smooth loop.  The
associated map \(S^1 \to L^x \R\) factors as \(S^1 \to L \R \to L^x
\R\).  The first factor is smooth by assumption whilst the second is a
continuous linear map and hence smooth.
\end{proof}

Thus although for continuity there is not much to choose between
\(L \R\) and \(L^0 \R\), once we get to smoothness we
easily see the difference.

\end{document}